\documentclass[a4paper,12pt,reqno]{amsart}

\usepackage[letterpaper, hmargin=1in, top=1in, bottom=1.2in, footskip=0.7in]{geometry}

\usepackage{rotating,pdflscape,xcolor,amssymb,subfigure,psfrag,amsmath,eufrak,bbm,epsfig,amsthm,mathtools,fancyhdr,graphicx}
\definecolor{vdarkred}{rgb}{0.6,0,0.2}
\definecolor{vdarkblue}{rgb}{0,0.2,0.6}
\usepackage[pdftex, colorlinks, linkcolor=vdarkblue,citecolor=vdarkred]{hyperref}
\usepackage{a4wide,amscd}
\usepackage{tikz} 
 \usepackage[usenames,dvipsnames]{pstricks}
 \usepackage{pst-grad} 
 \usepackage{pst-plot} 

\usepackage[small]{caption}

\newcommand{\ii}{\mathrm{i}}
\newcommand{\me}{\mathrm{e}}

\newcommand{\cB}{\mathcal{B}}

\newcommand{\cH}{\mathcal{H}}

\newcommand{\cM}{\mathcal{M}}
\newcommand{\cO}{\mathcal{O}}

\newcommand{\cT}{\mathcal{T}}
\newcommand{\cU}{\mathcal{U}}

\newcommand{\bv}{\mathbf{v}}
\newcommand{\bw}{\mathbf{w}}

\newcommand{\wtt}{\widetilde{t}}

\newcommand{\wTta}{\widetilde{\Tta}}

\newcommand{\si}{\sigma}
\newcommand{\lam}{\lambda}

\newcommand{\vfi}{\varphi}
\newcommand{\al}{\alpha}
\newcommand{\tta}{\theta}
\newcommand{\Tta}{\Theta}

\newcommand{\del}{\delta}
\newcommand{\Del}{\Delta}

\newcommand{\ld}{\ldots}
\newcommand{\cdote}{\,\cdot\,}

\newcommand{\beg}{\begin}
\newcommand{\en}{\end}

\renewcommand{\Re}{\mathfrak{Re}}

\newcommand{\trm}{\textrm}

\newcommand{\bgt}{\begin{itemize}}
\newcommand{\ent}{\end{itemize}}

\newcommand{\eqre}{\eqref}
\newcommand{\re}{\ref}
\newcommand{\la}{\label}
\newcommand{\ds}{\displaystyle}

\newcommand{\brem}{\begin{rmk}}
\newcommand{\erem}{\end{rmk}}
\newcommand{\blem}{\begin{lem}}
\newcommand{\elem}{\end{lem}}
\newcommand{\bcor}{\begin{cor}}
\newcommand{\ecor}{\end{cor}}
\newcommand{\bTh}{\begin{Th}}
\newcommand{\eTh}{\end{Th}}
\newcommand{\bpropo}{\begin{propo}}
\newcommand{\epropo}{\end{propo}}

\newcommand{\lan}{\langle}
\newcommand{\ran}{\rangle}

\newcommand{\op}{\operatorname}

\newcommand{\Tr}{\operatorname{Tr}}

\newcommand{\ud}{\mathrm{d}}

\newcommand{\ensn}{\{1,\ldots,n\}}

\newcommand{\E}{\op{\mathbb{E}}}

\newcommand{\R}{\mathbb{R}}
\newcommand{\C}{\mathbb{C}}

\newcommand{\p}{\mathbb{P}}

\newcommand{\pro}{probability }

\newcommand{\f}{\frac}
\newcommand{\ff}{\frac{1}}
\newcommand{\lf}{\left}
\newcommand{\ri}{\right}

\newcommand{\st}{such that }

\newcommand{\ti}{\times}

\newcommand{\ste}{\, ;\, }

\newcommand{\eps}{\varepsilon}

\newcommand{\bck}{\backslash}

\newcommand{\ovl}{\overline}

\newcommand{\bbm}{\begin{bmatrix}}
\newcommand{\ebm}{\end{bmatrix}}
\newcommand{\bes}{\begin{equation*}}
\newcommand{\ees}{\end{equation*}}
\newcommand{\be}{\begin{equation}}
\newcommand{\ee}{\end{equation}}
\newcommand{\beqy}{\begin{eqnarray}}
\newcommand{\eeqy}{\end{eqnarray}}
\newcommand{\beq}{\begin{eqnarray*}}
\newcommand{\eeq}{\end{eqnarray*}}
\newcommand{\one}{\mathbbm{1}}

\newcommand{\ie}{i.e. }

\newcommand{\bpm}{\begin{pmatrix}}
\newcommand{\epm}{\end{pmatrix}}

\newcommand{\cd}{\cdots}

\newcommand{\bpr}{\beg{proof}}
\newcommand{\epr}{\en{proof}}

\newcommand{\AND}{\qquad\trm{ and }\qquad}

 \newcommand{\theo}[1]{Theorem \re{#1}}

\newcommand{\lemm}[1]{Lemma \re{#1}}

\newcommand{\rem}[1]{Remark \re{#1}}

\newtheorem{Th}{Theorem}[section]

\newtheorem{propo}[Th]{Proposition}

\newtheorem{lem}[Th]{Lemma}

\newtheorem{cor}[Th]{Corollary}

\theoremstyle{definition}
\newtheorem{rmk}[Th]{Remark}
\newtheorem{conj}[Th]{Conjecture}

\long\def\symbolfootnote[#1]#2{\begingroup
\def\thefootnote{\fnsymbol{footnote}}\footnote[#1]{#2}\endgroup} 

\author{Florent Benaych-Georges} \address{Florent Benaych-Georges: MAP 5, UMR CNRS 8145 - Universit\'e Paris Descartes, 45 rue des Saints-P\`eres 75270 Paris cedex~6,  France.} \email{florent.benaych@gmail.com}
 \author{Ofer Zeitouni}
\address{Ofer Zeitouni: Department of Mathematics, Weizmann Institute of Science
  POB 26, Rehovot 76100, Israel
  and
  Courant Institute, New York University
 251 Mercer St, New York, NY 10012, USA. Supported in part by the ERC Advanced grant LogCorrelated-Fields.} \email{ofer.zeitouni@weizmann.ac.il}



\date{\today}
\title[Eigenvectors of non normal random matrices]{Eigenvectors of non normal random matrices}
\begin{document}
\maketitle
\beg{abstract}We study the angles between the eigenvectors of a random $n\times n$ complex matrix $M$ with density $\propto \mathrm{e}^{-n\operatorname{Tr}V(M^*M)}$ and $x\mapsto V(x^2)$ convex. We prove that for unit eigenvectors $\mathbf{v},\mathbf{v}'$ associated with distinct eigenvalues $\lambda,\lambda'$ that are the closest to specified points $z,z'$ in the complex plane, the rescaled inner product $$\sqrt{n}(\lambda'-\lambda)\langle\mathbf{v},\mathbf{v}'\rangle$$ is uniformly sub-Gaussian, and give a more precise statement in the case of the Ginibre ensemble.\en{abstract}

 \section{Introduction and main results}

 \subsection{Setup and main results}
 Let $V:\R^+\to \R$ be a function \st the following holds.
\begin{equation}
  \label{eq-alpha}
  \mbox{\rm For some $\al>0$, the function $V(x^2)-\f\al2{x^2}$ is convex.}
\end{equation}
 Let $X$ be an $n\ti n$ complex matrix with law  
 \begin{equation}
   \label{eq-law}
   \propto \me^{-n\Tr V(M^*M)}\ud M,
 \end{equation} 
 where 
 $\ud M$ is the standard Lebesgue measure on $n\ti n$ complex matrices. In particular,
 all eigenvalues of $X$ are distinct, almost surely (see \rem{whyeigenvalyesaredistinct}).
 
 Let $z,z'\in \C$ and let $\lam$ and $\lam'$ denote the eigenvalues of $X$ that are the closest to respectively  $z$ and $z'$ (if $z=z'$ or $\lam$ is the closest eigenvalue to both $z$ and $z'$, then $\lam'$ 
 is the second closest to $z'$).   
 Let  $\bv$ and $\bv'$ denote some associated eigenvectors of unit $\ell_2$ norm.  We want to study the quantity $n|\lam'-\lam||\lan \bv,\bv'\ran|^2$, which leads us to 
introduce the random variable $Y$, defined through any of the two following equivalent equations  \be\la{def:Yt12} Y\;:=\;  n \f{|\lam'-\lam|^2}{|\lan\bv,\bv'\ran|^{-2}-1}\ee  or 
\be\la{def:Yt12_2} n|\lam'-\lam|^2 |\lan\bv,\bv'\ran|^2\;=\;\f{ Y}{ \f{Y}{n|\lam'-\lam|^2}+1}, \qquad Y>0\ee 
Since eigenvalues are almost surely distinct and
$|\lan \bv, \bv'\ran|$ is invariant under multiplication of the eigenvectors by 
a complex scalar of norm $1$,  the random variable $Y$ is well defined.

Recall $\alpha$
from \eqref{eq-alpha}.
Our first main result is the following.

\beg{Th}\la{thCVeigvectors}The random variable $Y$ satisfies
\be\la{eqYt12}\p(Y\ge\del)\;\le\;2 \exp\lf(-\f{\al}2\del\ri), \quad
\mbox{\rm for any $\del>0$} .\ee
\en{Th}
In particular,
when $n|\lam'-\lam|^2\gg 1$, \eqre{def:Yt12_2} shows that
$$  |\lan\bv,\bv'\ran|\approx \f{\sqrt Y}{\sqrt n|\lam'-\lam|},$$
and it follows from Theorem \ref{thCVeigvectors} and a union bound that
all eigenvectors corresponding to mesoscopically separated
eigenvalues are asymptotically orthogonal to each other.

\begin{rmk}\label{rem-exponential}
In the case of the Ginibre ensemble ($V(x)=x$, so that the entries of 
$X$ are i.i.d. standard complex Gaussian variables with variance $n^{-1}$), 
the random variable
$Y$ has an exponential law of mean $1$. This fact is probably well known, 
and  
follows from Equations  \eqre{eq:fromYtot12}  and 
\eqre{densityTproof}-\eqre{def:ttilde} below. In particular,
$\sqrt{Y}$ is distributed like the norm of a standard complex Gaussian variable.
\end{rmk}
(We recall that  a \emph{standard complex Gaussian variable} is a  
centered complex Gaussian random variable  $Z$ \st $\E Z^2=0$ and $\E |Z|^2=1$.)

Our second main result is concerned with Ginibre matrices, for which
we extend an asymptotic version of Remark
\ref{rem-exponential} to the multivariate framework. 

\beg{Th}\la{thCVeigvectorsGinibre}Suppose that $V(x)=x$. For a  fixed $k\ge 2$,  let $z_1, \ld, z_k$ be (deterministic) points
in the unit disk, possibly dependent on $n$, 
\st for a certain $\eps>0$, uniformly in $n$, \be\la{300518111}  \sqrt n\min_{1\le i<j\le k} |z_j-z_i|\ge n^\eps .\ee
For each $i$, let $\lam_i$ be the eigenvalue that is the closest to $z_i$ and
let $\bv_i$ be 
an associated eigenvector. Let $\theta_i$, $i=1,\ldots,k$, be i.i.d. 
variables uniformly distributed on the $[0,2\pi]$, independent of $X$.
Then the distribution of the     triangular array 
$$\lf(\sqrt n(\lam_j-\lam_i) \lan \me^{\ii\theta_i} \bv_i,
\me^{\ii\theta_j}\bv_j\ran\ri)_{1\le i<j\le k}$$  
converges, as $n\to\infty$,  to the distribution of a triangular array of independent standard complex Gaussian variables.
\en{Th}

 \brem The typical distance between two eigenvalues of $X$ that are ``neighbors" of each other in the spectrum of $X$ has order $n^{-1/2}$. Hence, because of Hypothesis \eqre{300518111}, this result is well adapted for most pairs of eigenvalues, but not for those that are as close as possible  (in our proof, Hypothesis \eqre{300518111} is necessary for  estimates \eqre{eq:maxinv1805} to \eqre{300518c} to hold). For one given   pair of eigenvalues at distance $\asymp n^{-1/2}$, much information is contained in the fact that the random variable $Y$ introduced in \eqre{def:Yt12} and \eqre{def:Yt12_2} has exponential distribution with mean one. If one considers not only one such pair, but an arbitrary finite number $k$ of eigenvalues that are at distances of order $n^{-1/2}$, the problem is less simple.
 \erem

 \subsection{Background} The study of eigenvectors of random Ginibre matrices
 seems to have been initiated in \cite{CM1}. For a matrix $X$,
 let $v_i$ (respectively, $w_i$) denote the left (respectively, right) 
eigenvectors
corresponding to eigenvalues $\lambda_i$, where
the normalization 
$\lan \bv_i,\bw_j\ran=\delta_{ij}$ is imposed. Using the Schur representation
 $X=UTU^*$ with $T$ upper triangular and $U$ unitary, 
 they
 computed, for the Ginibre ensemble,
 the correlations of eigenvectors and cross correlations of right
 and left eigenvectors, with special emphasis on the correlator
 \begin{equation}
   \label{eq-cor}
   O_{1,2}=\E (\lan \bv_1,\bv_2\ran
\lan \bw_1,\bw_2\ran).
\end{equation}
 Using the joint density of entries of $T$, the 
 evaluation of the latter correlations reduce to the evaluation of certain 
 Green functions. This point of view was recently significantly expanded to more general models in \cite{NowakTarnowski} (using diagrammatic methods), as well as 
 in \cite{CR}, where multi-points correlations are evaluated and related to
 two point correlations. We refer the reader to the introduction of
 \cite{CR} for further details and an extensive bibliography.
Recent works \cite{Fyo,BD} study distributional limits for condition numbers, as well as refined estimates for overlaps in the microscopic and mesoscopic regime.

Another very relevant recent work  is \cite{BNST}, which deals with matrices
$X$ with joint density of entries of the form \eqref{eq-law} (whithout
assuming the convexity of $V$). In this general setup,
the correlator $O_{12}$ from \eqref{eq-cor}
is computed.

 Our results, as well as \cite{CM,CR},
 build 
upon the evaluation of the joint distribution of the entries of $T$, see
\cite{Ginibre, TAO2, Mehta}; these derivations 
do not address explicitely the ordering of the diagonal 
elements in $T$; in our approach, we choose the
ordering as function of the full set of diagonal elements.
For this reason, we provide explicitly a proof 
of the joint distribution of entries.

Finally, we mention that general delocalization results for eigenvectors
of random non-Hermitian matrices with independent entries 
appear in \cite{RV}.

 \subsection{A conjecture}
The inner products appearing in Theorem \ref{thCVeigvectorsGinibre}
can be written in terms of the off diagonal 
entries of the upper triangular matrix $T$ 
in the Schur decomposition of $X$, see the proof of Theorem
\ref{thCVeigvectors} below. In the Ginibre case, these
entries (scaled by $\sqrt{n}$) 
are iid standard complex Gaussians, and for general $V$,  
it is still the case that
$Y=n|t_{12}|^2$, see \eqref{eq:fromYtot12} below. This leads us to the following.
\begin{conj}
  Under the assumptions of Theorem
  \ref{thCVeigvectors},  the sequence 
  $Y/\E Y$ converges in distribution,
  as $n\to\infty$, to the exponential law of parameter one.
  \end{conj}
 Some preliminary computations make the conjecture plausible. In addition,
 the simulations in
 Figure \re{Fig:Spectral_Resolution} are in agreement with the conjecture.
\begin{figure}[!ht]
\centering
\subfigure[$V(x)=x$ (Ginibre)]{
\includegraphics[scale=.35]{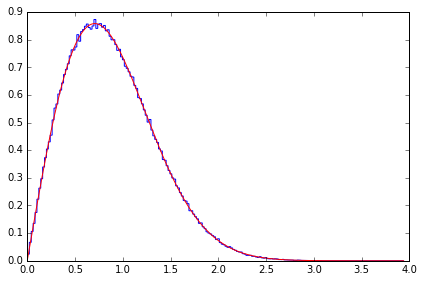}} \qquad 
\subfigure[$V(x)=x+\f{x^4}4+\f{x^5}5$]
{\includegraphics[scale=.35]{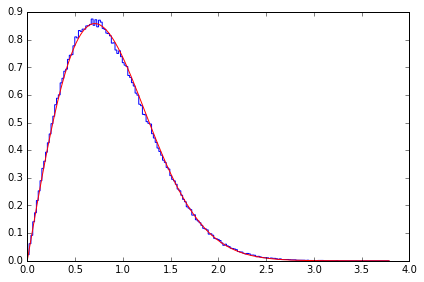}}
\caption{{\bf Universality of eigenvector angles:}  density   $\ds 2x\me^{-x^2}$ (in red) vs the histograms (in blue) of   $\ds   \f{\sqrt n|\lam'-\lam|}{\sqrt{|\lan\bv,\bv'\ran|^{-2}-1}}$ for $\lam\ne  \lam'$ running through the spectrum of a matrix with   distribution $\propto \me^{-n\Tr V(M^*M)}$ ($n=150$, sample of size $40$, sampled thanks to  Langevin Monte Carlo), with 2 different choices for $V$.  The sample on the right has been rescaled so that its empirical second moment is $1$ (because the distribution of $  \f{|\lam'-\lam|}{\sqrt{|\lan\bv,\bv'\ran|^{-2}-1}}$ can only be universal up to a rescaling of the matrix).
}\la{Fig:Spectral_Resolution}
\end{figure}

 \section{Proofs}
 In the proofs below, we use the joint distributions derived in
 Theorem \ref{thEigenvectors} from the appendix. 
 \bpr[Proof of \theo{thCVeigvectors}]  Set $\cO:=\cO_{z,z',z',z',\ld, z'}$ (see \eqre{defO} for the definition of this   set).  By Theorem \re{thEigenvectors}, we know that $X$ can be written $X=UTU^*$ with $U$ unitary and $T=[t_{ij}]$ upper triangular having the density  
 \be\la{densityT} \propto\one_{(t_{11},\ld, t_{nn})\in \cO }|\Del(t_{11},\ld,t_{nn})|^2 \me^{-n\Tr V(T^*T)},\ee
Hence by definition of 
$\cO$, $\lam$ and $\lam'$ are the two first diagonal entries of $T$. 
Thus the vectors $\bw:=(1, 0, \ld, 0)$ and $\bw':=\lf(|t_{12}|^2+|\lam'-\lam|^2\ri)^{-1/2}(t_{12}, \lam'-\lam, 0, \ld, 0)$ 
 are  unit right eigenvectors of $T$ for the eigenvalues $\lam$ and $\lam'$. 
Hence as $U$ is unitary, \be\la{CVproblamz0}|\lan\bv,\bv'\ran|=|\lan\bw,\bw'\ran|=\f{|t_{12}|}{\lf(|t_{12}|^2+|\lam'-\lam|^2\ri)^{1/2}},\ee so that 
\be\la{eq:fromYtot12} Y=   n \f{|\lam'-\lam|^2}{{|\lan\bv,\bv'\ran|^{-2}-1}}=  n |t_{12}|^2\ee
Thus we have to prove that   for any $\del>0$,  \be\la{eqYt12proof0}\p(\sqrt n|t_{12}|\ge\del)\;\le\;2\me^{-\f{\al\del^2}2}  .\ee 

Notice first that for any fixed  $(t_{11}, \ld, t_{nn})\in \cO$, the distribution $\mu_{\cdote|(t_{11}, \ld, t_{nn})}$  of $(t_{ij})_{1\le i<j\le n}$  conditional to $(t_{11}, \ld, t_{nn})$ has  on $\C^{n(n-1)/2}$ a density  \be\la{densityTconditionaldiag}\propto \me^{-n\Tr V(T^*T)}\ee for $T$ the upper-triangular matrix with upper-triangular entries $(t_{ij})_{1\le i\le j\le n}$. Thus by \lemm{lem:2505182} and \rem{rem:2505182}, $\mu_{\cdote|(t_{11}, \ld, t_{nn})}$ satisfies a LSI with constant $(\al n)^{-1}$.
Note also  that for any fixed $\tta_1, \ld,\tta_n\in \R$, the density of $\mu_{\cdote|(t_{11}, \ld, t_{nn})}$ given at \eqre{densityTconditionaldiag} is invariant under the transformation $$(t_{ij})_{1\le i<j\le n}\longmapsto (\me^{\ii(\tta_i-\tta_j)}t_{ij})_{1\le i<j\le n}.$$ We deduce that the expectation of $t_{12}$ with respect to $\mu_{\cdote|(t_{11}, \ld, t_{nn})}$ vanishes, 
so that by the LSI,   for any $\del>0$,  \be\la{eqYt12proof}\mu_{\cdote|(t_{11}, \ld, t_{nn})}(\sqrt n|t_{12}|\ge\del)\;\le\;2\me^{-\f{\al\del^2}2}  .\ee Integrating over $(t_{11}, \ld, t_{nn})$, we get \eqre{eqYt12proof0}. 
 \epr
 
 Before the proof of \theo{thCVeigvectorsGinibre}, we
 prove two preliminary lemmas. We suppose here that $V(x)=x$.

  \blem\la{lem:spacing_eigprelim} Let $z_0$ in the unit disk, possibly depending on $n$ and $s\in (0, 1/2)$. Then with \pro tending to one  as $n\to\infty$, $X$ has at least $\f{n^{1-2s}}5$ eigenvalues at distance $\le n^{-s}$ from $z_0$.
  \elem
  
  \bpr   Let $N_s$ denote the number of eigenvalues in the disk $D(z_0, n^{-s})$.  Let $f$ be a smooth non negative function    with value $1$ on the   disk $D(0,1/2)$ and with support contained in the disk $D(0,1)$. 
  Then we have $$N_s\ge \sum_{j=1}^nf(n^s(\lam_j-z_0)),$$ where the $\lam_j$ denote the eigenvalues.
  By the local circular law by Yin \cite[Th. 1.2]{yinLCL3} (see also \cite[Th. 9]{TaoVu} for the case where $|z_0|<1$), we know that with \pro tending to one, 
  $$\sum_{j=1}^nf(n^s(\lam_j-z_0))\ge \f n\pi \int_{|z|\le 1}f(n^s(\lam_j-z_0))\ud L(z)-n^{\si},$$ where $L$ denotes the Lebesgue measure on $\C$ and $\si:=(1-2s)/2$. 
  We deduce that with \pro tending to one, $$N_s\ge\sum_{j=1}^nf(n^s(\lam_j-z_0))\ge \f {n^{1-2s}}4-n^\si.$$
  \epr
  
     \blem\la{lem:spacing_eig}  Let $z,z'$ in the unit disk, possibly depending on $n$, \st for a certain fixed $\eps\in (0, 1/2)$, uniformly in $n$, $$\sqrt n |z'-z|\ge n^\eps. $$
Let $\lam, \lam'$ be the eigenvalues of $X$ that are the    closest to respectively  $z$ and $z'$ (if  $\lam$ is the closest eigenvalue to both $z$ and $z'$, then $\lam'$ is the second closest to $z'$). Then for any fixed $\del\in (0, \eps)$, we have  $$\sqrt n |\lam'-\lam|\ge n^\del$$ with \pro tending to one  as $n\to\infty$. \elem

\bpr Let  $s\in (1/2-\eps,1/2)$. By the previous lemma, we have $$|\lam-z|\le n^{-s}\AND |\lam'-z'|\le n^{-s}$$ with \pro tending to one as $n\to\infty$.
Thus $$|\lam'-\lam|\ge |z'-z|-|\lam-z|-|\lam'-z'|\ge n^{-1/2+\eps}-2n^{-s},$$ which allows to conclude, as $-s<-1/2+\eps$ and $ \eps>\del$.
\epr
 
 \brem In Lemmas \re{lem:spacing_eigprelim} and \re{lem:spacing_eig}, the properties do not only hold with \pro tending to one but  with \pro at least $1-Cn^{-D}$
 for any $D>0$ (and for $C$ a constant depending only on $D$, not on $z_0, z,z'$), which can be useful when using a union bound. The proof is the same and follows from the fact that in \cite{yinLCL3}, the error \pro is $\le Cn^{-D}$. \erem
 
  \bpr[Proof of \theo{thCVeigvectorsGinibre}] 
  For $u,v$ some random variables implicitly depending on $n$, we use the notation $u\sim v$ (resp. $u=O(v)$, $u=o(v)$) when $u/v$ tends in \pro to one (resp.  $u/v$ is tight,  $u/v$ tends in \pro to $0$) as $n\to\infty$. 

 Set $\cO:=\cO_{z_1,\ld,z_k,z_k, \ld, z_k}$ (see \eqre{defO} for the definition of this   set).  By Theorem \re{thEigenvectors}, we know that $X$ can be written $X=UTU^*$ with $U$ unitary and $T=[t_{ij}]$ upper triangular having the density  
 \be\la{densityTproof} \propto\one_{(t_{11},\ld, t_{nn})\in \cO }|\Del(t_{11},\ld,t_{nn})|^2 \me^{-n\sum_{1\le i\le j\le n}|t_{ij}|^2},\ee so that the random variables \be\la{def:ttilde}\wtt_{ij}:=\sqrt nt_{ij}\ee are  independent  standard complex Gaussian variables.

By definition of 
$\cO$, $\lam_1,\ld,\lam_k$ are the $k$ first diagonal entries of $T$.

Besides,  as $U$ is unitary, 
\be\la{eq:veqw}\lf(\sqrt n(\lam_j-\lam_i) \lan\bv_i,\bv_j\ran\ri)_{1\le i<j\le k}=\lf(\sqrt n(\lam_j-\lam_i) \lan\bw_i,\bw_j\ran\ri)_{1\le i<j\le k}\ee where the $\bw_i$ are the eigenvectors of $T$ associated to the $\lam_i$ (multiplied by independent uniform phases $\me^{\ii\tta_i}$, independent of $T$).

For each $i$, $\bw_i$ is in the kernel of $T-\lam_i$, hence has only its $i$ first coordinates non zero, and these coordinates are proportional to the vector $(x_i(1), \ld, x_i(i))\in \C^i$, satisfying \beq
(\lam_{1}-\lam_i)x_i(1)+t_{1,2}x_i(2)+\cdots\cdots\cdots\cdots\cdots\cdots\cdots+t_{1,i}x_i(i)&=&0\\ 
\ddots\qquad \qquad\qquad\qquad\qquad\qquad\qquad \qquad\qquad \qquad\vdots\;\;&=&\,\vdots\\
(\lam_{i-2}-\lam_i)x_i(i-2)+t_{i-2,i-1}x_i(i-1)+t_{i-2,i}x_i(i)&=&0\\
(\lam_{i-1}-\lam_i)x_i(i-1)+t_{i-1,i}x_i(i)&=&0\\
x_i(i)&=& 1\eeq
We solve this linear system: 
\beq &&x_i(i)\;=\; 1
\\
&&x_i(i-1)\;=\;\ff{\lam_i-\lam_{i-1}}t_{i-1,i}\\
&&x_i(i-2)\;=\;\ff{\lam_i-\lam_{i-2}}\lf(t_{i-2,i-1}x_i(i-1)+t_{i-2,i}\ri)\\
&&x_i(i-3)\;=\;\ff{\lam_i-\lam_{i-3}}\lf(t_{i-3,i-2}x_i(i-2)+t_{i-3,i-1}x_i(i-1)+t_{i-3,i}\ri)\\
&&\;\vdots\\ &&\;\vdots\\
&&x_i(1)\;=\;\ff{\lam_i-\lam_{1}}\lf(t_{1,2}x_i(2)+\cdots+t_{1,i-1}x_{i}(i-1)+t_{1,i}\ri)
\eeq
For each $1\le i<j\le k$, we have \be\la{eq:wxi}\lan \bw_i,\bw_j\ran=\f{\me^{\ii(\tta_j-\tta_i)}}{\|x_i\|\|x_j\|}\sum_{\ell=1}^i\ovl{x_i(\ell)}x_j(\ell).
\ee
To analyse the asymptotic behavior of these inner products, let us analyse the asymptotic behavior of each variable $x_i(\ell)$, $1\le \ell\le i\le k$. 

By \eqre{300518111} and  \lemm{lem:spacing_eig}, using the fact that $k$ is fixed, we have  
\be\la{eq:maxinv1805}\max_{1\le i<j\le k}\ff{|\lam_j-\lam_i|}=o(n^{\ff 2-\f\eps2}).\ee

By the previous equations and the estimate \eqre{eq:maxinv1805}, using the fact that the random variables  $\wtt_{ij}=\sqrt nt_{ij} $  are independent standard complex Gaussian variables, we have, for any $i=1, \ld, k$, we obtain successively the following estimates:
\beqy\nonumber &&x_i(i)\;=\; 1
\\
\la{300518a}&&x_i(i-1)\;\sim\;\ff{\lam_i-\lam_{i-1}}\f{\wtt_{i-1,i}}{\sqrt n}\;=\;o(n^{-\eps/2})\\
\la{300518b}&&x_i(i-2)\;\sim\;\ff{\lam_i-\lam_{i-2}}\f{\wtt_{i-2,i}}{\sqrt n}\;=\;o(n^{-\eps/2})\\
\nonumber&&\;\vdots\\ \nonumber&&\;\vdots\\
\la{300518c}&&x_i(1)\;\sim\;\ff{\lam_i-\lam_{1}}\f{\wtt_{1,i}}{\sqrt n}\;=\;o(n^{-\eps/2})
\eeqy
It implies that for all $1\le\ell\le i-1$, $$x_i(\ell)=O\lf(\ff{\sqrt n(\lam_i-\lam_\ell)}\ri)$$ and that, as $x_i(i)=1$,  
 $$\|x_i\|\sim 1.$$ By \eqre{eq:wxi}, we deduce that 
 $$\sqrt n(\lam_j-\lam_i)\lan \bw_i,\bw_j\ran\sim {\me^{\ii(\tta_j-\tta_i)}} \wtt_{ij}+ \ff{\sqrt n}O\lf(\sum_{\ell=1}^{i-1} \f{\lam_j-\lam_i}{\ovl{\lam_i-\lam_\ell}(\lam_j-\lam_\ell)}\ri).$$
 For each $\ell=1, \ld, i-1$, 
 $$\f{\lam_j-\lam_i}{\ovl{\lam_i-\lam_\ell}(\lam_j-\lam_\ell)}=\ff{\ovl{\lam_i-\lam_\ell}}-\f{\lam_i-\lam_\ell}{\ovl{\lam_i-\lam_\ell}}\ff{\lam_j-\lam_\ell}=o(n^{\ff 2-\f\eps2}),$$ where we used \eqre{eq:maxinv1805}. 
 It follows that $$\sqrt n(\lam_j-\lam_i)\lan \bw_i,\bw_j\ran\sim {\me^{\ii(\tta_j-\tta_i)}} \wtt_{ij}+ o(n^{-\eps/2}),$$
  and, as $(\wtt_{ij})_{1\le i<j\le k}$ is a collection of independent standard complex Gaussian variables independent of the $\tta_i$'s,  the result is proved.
  \epr
 \section{Appendix}
 \subsection{Change of variables in the Schur decomposition}
  We endow the sets $\cM_n(\C)$ and $\cT_n(\C)$ of $n\ti n$   respectively complex matrices and upper-triangular  complex matrices with the Euclidian structures defined by \be\la{def:EuclidianStructure}X\cdot Y:=\Re\Tr(XY^*)\ee and let $\ud M$ (resp. $\ud T$) denote the associated Lebesgue measure on $\cM_n(\C)$ (resp. on $\cT_n(\C)$).   We also denote by $\cU_n$ the group  of unitary $n\ti n$ matrices and by $\ud U$ the Haar measure on $\cU_n$. 
 
 \brem\la{whyeigenvalyesaredistinct} It is useful to note that the set of matrices  in $ \cM_n(\C)$ with multiple eigenvalues is the set of matrices whose characteristic polynomial has null discriminant (see \cite[Def. A.10]{Alice-Greg-Ofer}), so that this set is a level set of a non constant polynomial function on  $ \cM_n(\C)$, hence has zero Lebesgue measure (the last fact can be checked by applying Fubini's theorem).    
 \erem

We begin by 
defining \emph{admissible  sets}, a notion which will allow us to order the eigenvalues of non Hermitian matrices in quite general ways. $S_n$ denotes the set of permutations of $\ensn$.
 
  \beg{Def}\la{defadmiset} An open set $\cO\subset\C^n$ is said to be \emph{admissible} if the sets $\si\cdot \cO:=\{(t_{\si(1)}, \ld, t_{\si(n)})\ste (t_1, \ld, t_n)\in \cO\}$, $\si\in S_n$, are pairwise disjoint and $\C^n\bck\cup_{\si\in S_n}\si\cdot \cO$ has null Lebesgue measure.\en{Def}
  
  An important example of  
admissible set is the following one. 
Fix $z_1, \ld, z_n\in \C$. Then the set   of $n$-tuples $(t_1, \ld, t_n)\in\C^n$ where for each $i$, the $i$-th entry $t_i$ is  strictly  closer to $z_i$ than all  the forthcoming ones, \ie the set 
\be\la{defO} \cO_{z_1, \ld, z_n}:=\{(t_1, \ld, t_n)\in\C^n\ste \forall   i<j, |z_i-t_i|<|z_i-t_j|\},\ee
   is admissible.

 \beg{Th}\la{thEigenvectors}Let $\rho$ be a non negative measurable function on $\cM_n(\C)$ \st for any $M\in \cM_n(\C)$ and any unitary matrix $U$, \be\la{rhoinv}\rho(M)=\rho(UMU^*).\ee Fix an admissible set $\cO\subset \C^n$. Then the measure $\rho(M)\ud M$ on $\cM_n(\C)$ is the push-forward, by the function $(U,T)\mapsto UTU^*$, of the measure $$C_n\ud U\otimes \lf(\one_{(t_{11},\ld, t_{nn})\in \cO}|\Del(t_{11},\ld,t_{nn})|^2\rho(T)\ud T\ri)$$ on $\cU_n\times \cT_n(\C)$, where $$\Del(t_{11},\ld,t_{nn}):=\prod_{1\le i<j\le n}(t_{jj}-t_{ii})$$ and $C_n$ is a constant depending only on $n$ (and not on $\rho$). 
 \en{Th}
 
 \brem Using the case of Ginibre matrices, one can compte $C_n$ : \be\la{Eq:Cn}C_n=\ff{\pi^{\f{3n^2-n}2}\prod_{1\le k\le n-1}k!}
\ee
 \erem
 
 \bpr[Proof of \theo{thEigenvectors}]   Some statements which are very close to \theo{thEigenvectors} are proved in various texts, as \cite{Ginibre,Mehta,Forrester,TAO2}.  However, firstly, these results are a bit less general and written in slightly different languages and, secondly and more importantly, they do not treat the question of the ordering the diagonal entries of $T$ (which is the cornerstone of our approach in this paper). For this reason, we provide a complete proof.

 \beg{lem}\la{lem_fondamental}Let $\cT, \cU, \cM$ be some open subsets of respectively $\R^p$, $\R^q$ and $\R^{p+q}$. Let $\vfi : \cT\ti \cU\to \cM$ be a smooth diffeomorphism with reciprocal denoted by $\Psi=(\Psi_1, \Psi_2)$. Let also $\rho$ be a non negative measurable function on $\cM$. Let $\ud t, \ud u, \ud m$ denote the Lebesgue measures on respectively $\cT, \cU, \cM$. Then the push-forward of the measure $\rho(m)\ud m$ on $\cM$ by the function $\Psi_1 : \cM\to \cT$ is $K(t)\ud t$, with $$K(t):=\int_{u\in \cU}\rho\circ \vfi(t,u)|J\vfi(t,u)|\ud u,$$ with $|J\vfi(t,u)|$ the Jacobian\footnote{What we call here the \emph{Jacobian} of a smooth  function between two Euclidian spaces with the same dimension is the absolute value of the determinant of the matrix of its derivative in any pair of orthogonal bases.} of $\vfi$.
 \en{lem}
 
 \bpr Let $f: \cT\to\R$ be a test function. Then \beq \int_{m\in \cM} f\circ \Psi_1(m)\rho(m)\ud m&=&\int_{(t, u)\in \cT\ti\cU} \underbrace{f\circ \Psi_1\circ\vfi(t,u)}_{f(t)}\rho\circ \vfi(t,u)|J\vfi(t,u)|\ud u\ud t\\ &=&\int_{t\in \cT} f(t)K(t)\ud t.\eeq
 \epr
 
  \beg{Def}\la{defmatsets}  Let $\cU'$ be the set of   unitary matrices whose entries are all non zero, whose diagonal entries are positive and whose principal minors are all invertible. \en{Def} 
 
 The following lemma can be found in \cite[Lemma 2.5.6]{Alice-Greg-Ofer}.

\beg{lem}\la{defTheta}
The map $$\Xi:\cU'\to \R^{  n(n-1)}$$ which maps $U$ to \be\la{1521216h39}\bbm 0& \f{u_{12}}{u_{11}} & \f{u_{13}}{u_{11}}&& \cd&\f{u_{1n}}{u_{11}}\\ \\
&0& \f{u_{23}}{u_{22}}&& \cd&\f{u_{2n}}{u_{22}}\\ \\ 
&&0&&&\\
&&&\ddots&&\\
&&&&0&\f{u_{n-1,n}}{u_{n-1,n-1}}\\
0&0&0&\cd&0&0
\ebm\ee
(without the zeros) is diffeomorphism from $\cU'$ onto a subset of  
 $ \R^{ n(n-1)/2}$ with closed null mass complementary. We denote its inverse by $\Theta$.\en{lem}

 Set  \be\la{def:tcO}\cT_{\cO}:=\{T=[t_{ij}]\in \cT_n(\C)\ste  (t_{11},\ld, t_{nn})\in \cO\}.\ee
 
\beg{lem}[Schur decomposition]\la{SchurDecLemma}For $ \cU'$ as 
in Definition \re{defmatsets}  and $\cT_{\cO}$  as in \eqre{def:tcO}, 
there is an open subset $\cM'\subset\cM_n(\C)$ with   null mass complementary \st any $M\in \cM'$  can be written in a unique way \be\la{SchurEq}M=UTU^*\ee with $T\in \cT_{\cO}$   and $U\in \cU'$. 
\en{lem}
Note that if $M=UTU^*$, then 
  \be\la{eq29junerho}\rho(M)=\rho(T).\ee 

 \beg{Def}\la{defmatsets3} For $\Tta$ as defined in Lemma \re{defTheta}, let
$\vfi : \cT_{\cO}\ti \R^{n(n-1)}\to \cM'$ 
be the diffeomorphism defined by $\vfi(T,x)=\Tta(x)T\Tta(x)^{*}$ and let $\Psi=(\Psi_1, \Psi_2)$ be its inverse.\en{Def}
By the unitary invariance of  \eqre{rhoinv}, the proof of the theorem reduces to the proof of the fact that 
 the push-forward, by $\Psi_1$, of the measure  $$\rho(M)\ud M$$ is$$
 C_n \one_{(t_{11},\ld, t_{nn})\in \cO}\Del(t_{11},\ld,t_{nn})^2\rho(T)\ud T$$ for $C_n$   a constant depending only on $n$. 
 By \eqre{eq29junerho} and Lemma \re{lem_fondamental}, this push-forward is the measure   $$\lf(\one_{(t_{11},\ld, t_{nn})\in \cO}\rho(T)\int_{x\in \R^{n(n-1)}}|J\vfi(T,x)|\ud x\ri)\ud T.$$

Then, the following lemma concludes the proof of  \theo{thEigenvectors}.
\beg{lem}\la{lem:20015181}On $\cT_{\cO}\ti \R^{n(n-1)}$, we have $$|J\vfi(T, x)|= g(x)\prod_{1\le i<j\le n}|t_{jj}-t_{ii}|^2,$$ with  $g(x)$ a measurable function of $x$.\en{lem}

\bpr 
Let    $\wTta: \cT_n(\C)\ti \R^{ n(n-1)}\to \cT_n(\C)\ti\cU'$ be defined by  $\wTta(T, x):=(T, \Tta(x))$ and $F: \cT_n(\C)\ti \cU'\to \cM_n(\C)$ be defined by $F(T, U):=UTU^*$ (note that $F$ is defined on a manifold and not on an open subset of an Euclidian space). We have $\vfi=F\circ \wTta$  on $\cT_{\cO}\ti \R^{  n(n-1)}$, so we have $$|J\vfi(T, x)|=|JF(T, \Tta(x))|\ti |J\Tta(x)|,$$ hence it suffices to prove that on $\cT_{\cO}\ti \cU'$, $$|JF(T ,U)|=f(U)\prod_{1\le i<j\le n}|t_{jj}-t_{ii}|^2,$$ where $f(\cdote)$ is a function of $U$.

Note that the tangent space of $\cU'$ at $U$ is  the space $$\op{Tangent}_U(\cU')=\{UH\ste H\in\ii\cH_0\}$$
for 
\be\la{defiH0}\ii\cH_0:=\{H=[h_{ij}]\in \cM_n(\C)\ste H^*=-H\trm{ and }\forall i,   h_{ii}=0\}.\ee Note also   that for all $T,R\in \cT_n(\C)$, $U\in \cU'$ and $H\in \op{Tangent}_U(\cU')$, we have 
\beq DF(T, U)(R,H)&=&HTU^*+URU^*+UTH^*\\
&=&U(U^*HT+R+TH^*U)U^*\\
&=&U(U^*HT-TU^*H+R)U^*
\eeq
 As the transformation of $\cM_n(\C)$ defined by  $K\mapsto UKU^*$ is orthogonal for any unitary $U$, $|JF(T,U)|$ is the absolute value of the determinant of the matrix, in an orthonormal basis, of the map $ \cT_n(\C) \ti \ii\cH_0\to \cM_n(\C)$ defined by $(R,K)\mapsto KT-TK+R$.
Using the fact that $\cM_n(\C)$ is the orthogonal sum of $\cT_n(\C)$ and of the space \be\la{defiTsl}\cT^{sl}:=\{M=[m_{ij}]\in \cM_n(\C)\ste \forall i\le j, m_{ij}=0\}\ee of strictly lower triangular $n\ti n$ matrices, it is easy to see that the determinant of this map is the one of the   map $\mathfrak{C}_T$ of Lemma \re{2112954} below, which  concludes the proof of \lemm{lem:20015181}. 
\epr

\beg{lem}\la{2112954}Let  $\cT^{sl}$ be as in \eqre{defiTsl}, $\ii\cH_0$ be as in \eqre{defiH0} and $\pi:\cM_n(\C)\to \cT^{sl}$ be the canonical  projection. Then for any  $T=[t_{ij}]\in \cT_n(\C)$,  the map $\mathfrak{C}_T:  \ii\cH_0\to \cT^{sl}$ defined by $\mathfrak{C}_T(M):=\pi(MT-TM)$,  has Jacobian $$\prod_{1\le i<j\le n}|t_{jj}-t_{ii}|^2,$$ all spaces being endowed with the Euclidian structure induced by \eqre{def:EuclidianStructure}.
 \en{lem}
 
 \bpr To prove this lemma, we shall first fix some orthonormal bases of $\ii\cH_0$ and $\cT^{sl}$, order them and then prove that the matrix of $\mathfrak{C}_T$ on these (conveniently ordered) bases is lower triangular  by $2\ti 2$ blocs with diagonal blocs having determinants $|t_{jj}-t_{ii}|^2$, $1\le j<i\le n$.
 
 Let  us denote the elementary matrices by $E_{ij}$ and let $\cB_{\ii\cH_0}$ be the family $$\ff{\sqrt{2}}(E_{ij}-E_{ji},\ii E_{ij}+\ii E_{ji})_{1\le j<i\le n}$$ and let $\cB_{\cT^{sl}}$ be the family $$(E_{ij},\ii E_{ij})_{1\le j<i\le n}.$$ These are of course orthonormal bases of $\ii\cH_0$ and $\cT^{sl}$. Order the set   $\ensn^2$ with the lexicographical order made out of the reverse natural order on the first component  and the natural order on the second one (for this order, for example, $(1,3)\le (1,4)$ and  $(2,1)\le (1,5)$) and use the induced order on  
 the bases $\cB_{\ii\cH_0}$ and $\cB_{\cT^{sl}}$.

 Let $M=[m_{ij}]\in {\ii\cH_0}$  and let $[v_{ij}]:=MT-TM$.
 We have \be\la{fortrisuplexic}v_{ij}=\sum_{k=1}^jm_{ik}t_{kj}-\sum_{k=i}^n t_{ik}m_{kj}=m_{ij}(t_{jj}-t_{ii})+\sum_{k=1}^{j-1}m_{ik}t_{kj}-\sum_{k=i+1}^n t_{ik}m_{kj}\ee
It follows that the matrix of 	$\mathfrak{C}_T$ on  the bases $\cB_{\ii\cH_0}$ and $\cB_{\cT^{sl}}$ (ordered as above) is  lower diagonal by $2\ti 2$ blocs, with $2\ti 2$ diagonal blocs the matrices of the linear maps $\C\to\C$ ($\C $ considered as a real vector space) $m\longmapsto (t_{jj}-t_{ii})\ti m$, $1\le j<i\le n$. The determinant of such a map is $|t_{jj}-t_{ii}|^2$, so the result follows.
 \epr \epr

\subsection{Klein's lemma and consequences}
The following lemma can be found in \cite[Lemma 4.4.12]{Alice-Greg-Ofer}.
\blem[Klein's lemma] For any $f:\R\to \R$ convex and $n\ge 1$, the function $M\mapsto \Tr f(M)$ is convex on the space of $n\ti n$ Hermitian matrices.\elem 
 
 For the previous lemma, we shall deduce the following one.
 \blem\la{lem:2505181} Let $f:\R^+\to\R$ be \st $g(x):=f(x^2)$ is convex  and $n\ge 1$. Then the function $X\mapsto \Tr f(X^*X)$ is convex on the space $\cM_n(\C)$ of $n\ti n$ complex matrices.
\elem
\brem\la{rem:2505181} As a direct consequence, for any fixed $(t_{11}, \ld, t_{nn})\in \C^n$, the function $$(t_{ij})_{1\le i<j\le n}\in \C^{n(n-1)/2}\mapsto \Tr f(T^*T), $$ for $T$ the upper-triangular matrix with   entries $(t_{ij})_{1\le i\le j\le n}$,  is convex.
\erem

\brem\la{rem:2505182} Suppose now, with the notation of the lemma, that for some $\al>0$,    $g(x)-\f\al2{x^2}$ is convex. Then, by the lemma, the  function $X\mapsto \Tr f(X^*X)-\f{\al}2\Tr X^*X$ is convex on the space $\cM_n(\C)$ of $n\ti n$ complex matrices. In the framework of Remark \re{rem:2505181}, as adding a constant to a convex function doesn't break convexity, it implies that the function $$(t_{ij})_{1\le i<j\le n}\in \C^{n(n-1)/2}\mapsto \Tr f(T^*T)-\f{\al}2\sum_{1\le i<j\le n}|t_{ij}|^2$$ is convex.
\erem

\bpr[Proof of \lemm{lem:2505181}] By Klein's lemma, the function $$X\in\cM_n(\C)\mapsto \Tr g \lf(\bpm 0&X\\ X^*&0\epm\ri)$$ is convex. Then, conclude noting that for $M:=\bpm 0&X\\ X^*&0\epm$,  $$\Tr f(X^*X)=\ff2\lf(\Tr f(X^*X)+\Tr f(XX^*)\ri)=\ff2\Tr f(M^2)=\ff2\Tr g(M).$$ 
\epr

\subsection{Logarithmic Sobolev Inequalities and concentration}
The following lemma, due to Bobkov, Ledoux and Herbst,  gives a sufficient condition for a \pro measure to satisfy a logarithmic Sobolev Inequality (\emph{LSI})  and states one of  its main consequences (see \cite[Sec. 2.3.2]{Alice-Greg-Ofer} for a definition of LSI and   a reference for the following lemma).  

\begin{lem}\la{lem:2505182}For any $V: \R^n \to \R\cup\{+\infty\}$  and $\al>0$ \st  $V(x)-\f\al2{\|x\|_2^2}$ is convex,  the \pro measure $P_{V,\R^n}\propto \me^{-V(x)}\ud x$ satisfies a LSI with constant $\al^{-1}$. This implies that
 for any $1$-Lipschitz function $f:\R^n\to \R$ and any $\del>0$, we have
$$P_{V,\R^n}\lf(|f(x)-\E_{P_{V,\R^n}}f|\ge \del\ri)\;\le\;2\me^{-\al\del^2/2}.$$\end{lem}

  \begin{thebibliography}{10}
 \bibitem{Alice-Greg-Ofer} G.~Anderson, A.~Guionnet, O.~Zeitouni \emph{An Introduction to Random Matrices}. Cambridge studies in advanced mathematics, {118} (2009).
 \bibitem{BESLocal} Z. Bao, L. Erdos and K. Schnelli \emph{Local single ring theorem on optimal scale}, Ann. Probab. (to appear).
 \bibitem{BNST} S. Belinschi, M. A. Nowak, R. Speicher,
   W. Tarnowski
   \emph{Squared eigenvalue condition numbers and eigenvectors from the single ring theorem},  J. Phys. A: Math. Theor. 50 (2017), 105204.
 \bibitem{BD} P. Bourgade and G. Dubach, \emph{The distribution of overlaps 
   between eigenvectors of Ginibre matrices}, arXiv:1801.01219, (2018).
 \bibitem{CR} N. Crawford and R. Rosenthal \emph{Eigenvector correlators in 
   the complex Ginibre Ensemble}, arXiv:1805.08993 (2018).
 \bibitem{Forrester}  P. J. Forrester \emph{Log-gases and random matrices}, London Mathematical Society Monographs Series, vol. 34, Princeton University Press, Princeton, NJ, 2010.
 \bibitem{Fyo} Yan V. Fyodorov \emph{On statistics of bi-orthogonal eigenvectors in real and complex Ginibre ensembles: combining partial Schur decomposition with supersymmetry}, 	arXiv:1710.04699 (2017).
\bibitem{Ginibre} J.  Ginibre \emph{Statistical ensembles of complex, quaternion, and real matrices}. J. Mathematical Phys. 6 1965 440--449.
  \bibitem{GKZ} A. Guionnet, M. Krishnapur, O. Zeitouni \emph{The Single Ring Theorem}. Ann. of Math. (2) 174 (2011), no. 2, 1189--1217.
 \bibitem{CM1} B. Mehlig and J. T. Chalker \emph{Eigenvector correlations
   in non-Hermitian random matrix ensembles}. Ann. Phys. 7 (1998), 427--436.
    \bibitem{CM} B. Mehlig and J. T. Chalker \emph{Statistical properties of
     eigenvectors in non-Hermitian Gaussian non Hermitian Gaussian
   random matrices ensembles}. J. Math. Phys.  41 (2000), 3233-3256.
 \bibitem{Mehta} M. Mehta \emph{Random matrices}. Third Ed., Academic Press (2004).
 \bibitem{NowakTarnowski} M.A. Nowak,  W. Tarnowski  \emph{Probing non-orthogonality of eigenvectors in non-Hermitian matrix models: diagrammatic approach}. J. High Energ. Phys. (2018) 2018: 152.
 \bibitem{RV} M. Rudelson and R. Vershynin
   \emph{Delocalization of eigenvectors of random matrices with independent
   entries}, Duke Math. J. 164 (2015), 2507--2538.
 \bibitem{TAO2} T. Tao 	\emph{Topics in random matrix theory}, Graduate Studies in Mathematics, AMS, 2012. 
 \bibitem{TaoVu} T. Tao, V. Vu 	\emph{Random matrices: Universality of local spectral statistics of non--Hermitian matrices}, Annals of Probability 2015, Vol. 43, No. 2, 782--874.
 \bibitem{yinLCL3} I. Yin \emph{The local circular law III: general case}. Probab. Theory Relat. Fields (2014) 160:679--732.
 \en{thebibliography}
\en{document}